\documentclass{article}
\usepackage{spconf,amsmath,graphicx}

\usepackage{dsfont}
\usepackage{amssymb}
\usepackage{amsbsy}
 \usepackage{algorithmic}
 \usepackage{algorithm}
 \usepackage{url}
\usepackage{scalefnt}

\newcommand{\algorithmicoutput}{\textbf{Output:}}
\newcommand{\OUTPUT}{\item[\algorithmicoutput]}

\newcommand{\algorithmicinput}{\textbf{Input:}}
\newcommand{\INPUT}{\item[\algorithmicinput]}

\def\x{{\mathbf x}}
\def\y{{\mathbf y}}
\def\z{{\mathbf z}}
\def\u{{\mathbf u}}

\def\h{{\mathbf h}}
\def\V{{\mathbf V}}

\def\R{{\mathbb R}}

\def\xib{{\boldsymbol{\xi}}}

\title{Texture Modeling by Gaussian fields \\with prescribed local orientation }
\name{K\'evin Polisano$^1$, Marianne Clausel$^1$, Val\'erie Perrier$^1$, Laurent Condat$^2$\thanks{The authors acknowledge the support of the French Agence Nationale de la Recherche (ANR) under reference
ANR-13-BS03-0002-01 (ASTRES).}}
\address{$^1$University of Grenoble-Alpes, Laboratoire Jean Kuntzmann, UMR 5224 CNRS, Grenoble, France
\\ $^2$University of Grenoble-Alpes, 
GIPSA-lab, UMR 5216 CNRS, Grenoble, France}
\begin{document}
%
\def\baselinestretch{0.987}
\maketitle
\scalefont{0.953}
%
\begin{abstract}
This paper presents a new framework for oriented texture modeling. 
We introduce a new class of Gaussian fields, called {\it Locally Anisotropic Fractional Brownian Fields}, with prescribed local orientation at any point. These fields are a local version of a specific class of anisotropic self-similar Gaussian fields with stationary increments. The simulation of such textures is obtained using a new algorithm mixing the tangent field formulation and a turning band method, this latter  method having proved its efficiency for generating stationary anisotropic textures.
Numerical experiments show the ability of the method for synthesis of textures with prescribed local orientation.
\end{abstract}
\section{Introduction}
\label{sec:intro}

Texture modeling is a challenging issue of image processing. There is a variety of texture methods in the field of computer vision, namely structural, statistical, model-based and transform-based methods. Thus, identifying the perceived characteristics of a texture in an image (regularity, roughness, frequency, content directionality, etc.) is an important first step towards building mathematical models for textures. We are interested in textures presenting same similar patterns at different scales, as is often the case for objects appearing in the nature, like clouds or mountains. We focus on stochastic models with a property of self-similarity, characteristic of a fractal behavior. The stochastic model behind fractal analysis is the fractional Brownian field (FBF), which is a multi-dimensional extension of the famous fractional Brownian motion (FBM) introduced in 1940 by Kolmogorov \cite{kol1940} as a way to generate Gaussian ``spirals'' in Hilbert spaces. The systematic study of the FBM started with the 
seminal paper of Mandelbrot and Van Ness \cite{mandelbrot1968}. The FBM has now become a standard model: it is used in many areas such as hydrology, economics, finance, physics and telecommunications (see, e.g., \cite{samorodnitsky1997}, \cite{LV1997}, \cite{abry2002}, and references therein for more details). The FBF has also been largely used in medical applications, with for instance the study of lesion detectability in mammogram textures \cite{caldwell1990}, assessment of breast cancer risk \cite{heine2002}, and the characterization of bone architecture for the evaluation of osteoporotic fracture risk \cite{benhamou2001}. 

Nevertheless, in many cases, fractal analysis with fractional Brownian fields, which are isotropic by definition, is not completely satisfactory, in particular when the considered data display some anisotropy. Therefore, many stochastic models have been introduced in the literature to take into account these possible additional anisotropic properties. Let us cite notably fractional Brownian sheet defined in~\cite{kamont1996,ayache2002} and anisotropic fractional Brownian field (AFBF) introduced by Bonami and Estrade in~\cite{Bonami2003} which are two classical examples of Gaussian fields satisfying {\it global} anisotropic properties. Other models of anisotropic textures called locally parallel textures, have also been recently introduced in \cite{maurel2011}. The mathematical definition and computational synthesis of anisotropic textures is an important issue, since it provides flexible models enabling to test estimation procedures of the anisotropic characteristics of an image. Here we focus on 
anisotropic {\it local properties} of Gaussian textures and provide a new Gaussian model whose anisotropic properties are prescribed at every point. It is a first preliminary and important step in defining new statistical estimators of the local anisotropic features of a given texture.


The paper is organized as follows. Section 2 briefly reviews definitions and characterizations of a class of self-similar Gaussian fields derived from the AFBF. Section 3 is devoted to the presentation of our model, from both the theoretical and implementation point of view. Finally, we provide the synthesis of numerical textures, for several vector fields of local orientations, showing the ability of our approach. 

\section{Anisotropic self-similar Gaussian fields}
\label{sec:afbf}

%


\subsection{The fractional Brownian field}
\begin{figure}[t]
\begin{minipage}[b]{.48\linewidth}
  \centering
  \centerline{\includegraphics[width=4.0cm]{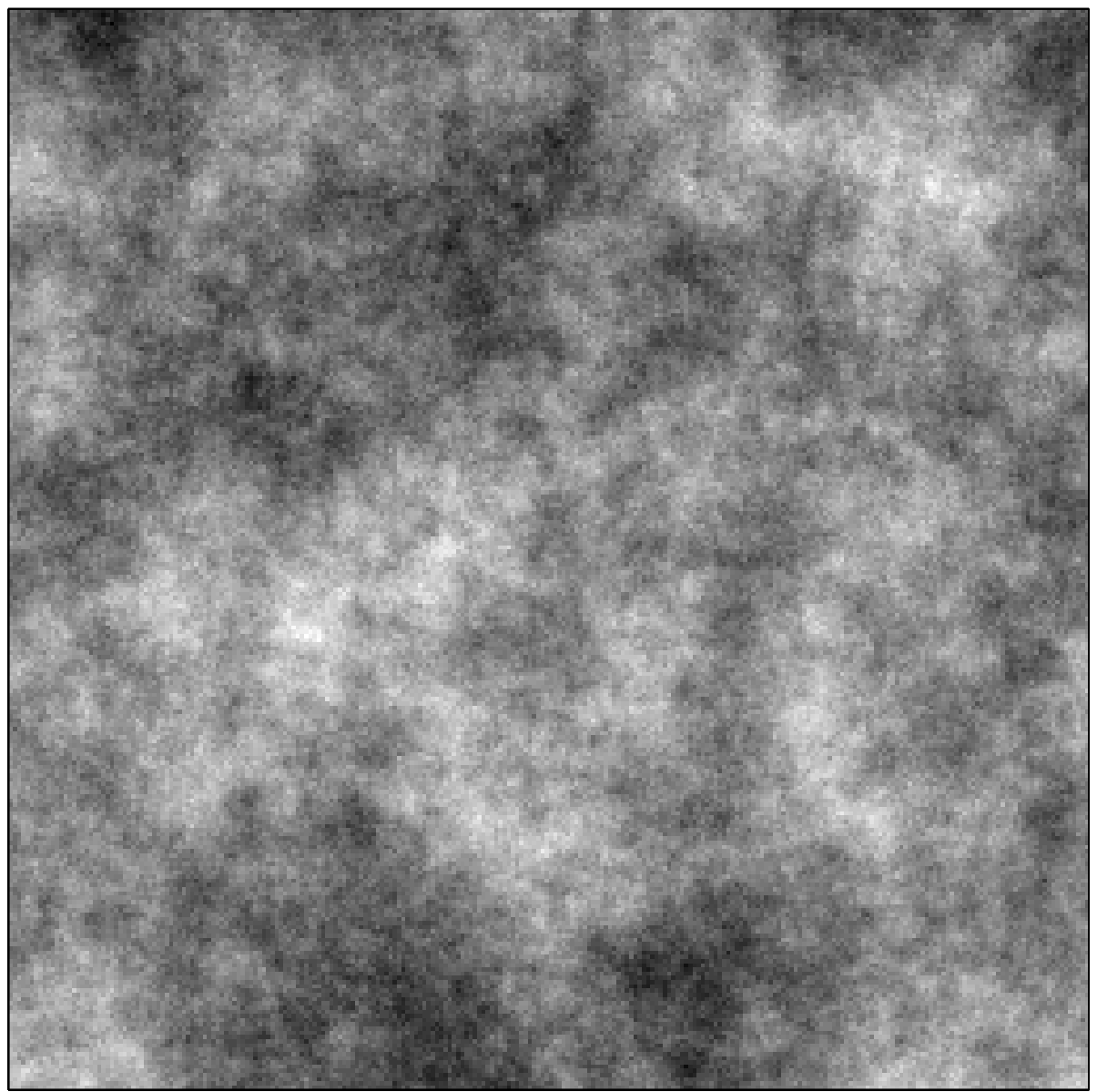}}
  \centerline{(a)}\medskip
\end{minipage}
\hfill
\begin{minipage}[b]{0.48\linewidth}
  \centering
  \centerline{\includegraphics[width=4.0cm]{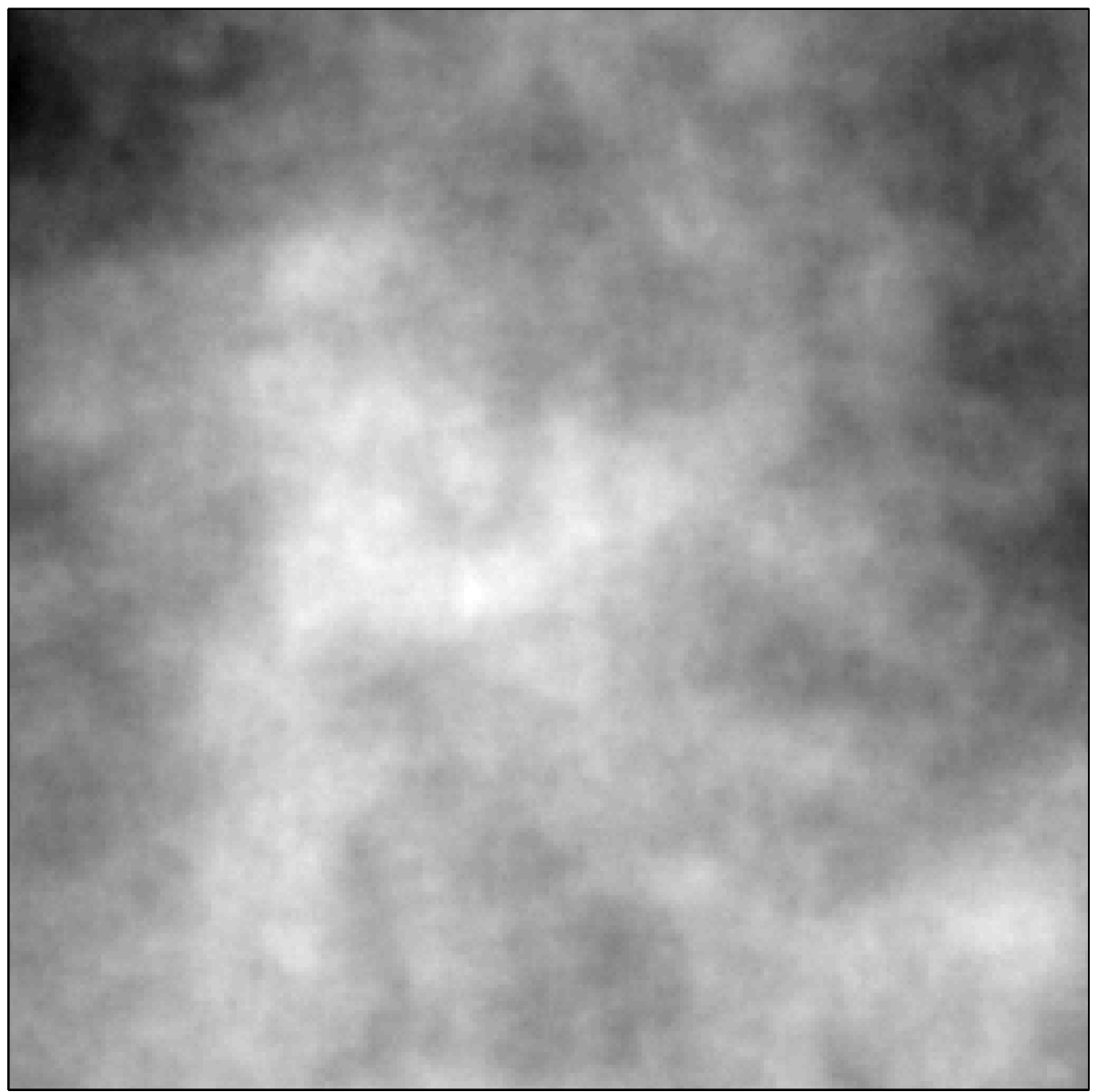}}
  \centerline{(b)}\medskip
\end{minipage}
\caption{Examples of FBF with (a) $H=0.3$, (b) $H=0.7$.}
\label{fig:fbf}
\end{figure}

Let $0<H<1$. The \textit{fractional Brownian field} with Hurst index $H$, denoted by
$B^H=\{B^H(\x);\x\in \R^2\}$, is the unique real-valued centered Gaussian field satisfying the following
properties:
 -- almost surely $B^H(0)=0$,\\
 -- $B^H$ admits stationary increments, i.e, for every $\z\in \R^2,$ $B^H(\cdot + \z)-B^H(\z)\overset{\mathcal{L}}{=} B^H(\cdot)-B^H(0),$\\
 -- $B^H$ is $H$ self-similar, i.e, $\forall \lambda \in \R^{\star}, B^H(\lambda \cdot)\overset{\mathcal{L}}{=} \lambda^{H}B^H(\cdot),$\\
 -- $B^H$ is isotropic, i.e, for every rotation $R$ in $\R^2$, $B^H\circ R\overset{\mathcal{L}}{=} B^H,$\\
where $\overset{\mathcal{L}}{=}$ denotes equality for all finite dimensional distributions. The FBM is wholly characterized by its covariance function, which is given, for every $\x,\y\in \R^2$ by $$\text{Cov}(B^H(\x),B^H(\y))=c_{H}(\|\x\|^{2H}+\|\y\|^{2H}-\|\x-\y\|^{2H})\;,$$
$c_H$ being a well-known nonnegative constant depending on $H$. Following \cite{samorodnitsky1997}, the FBM can also be defined by its harmonizable representation:
\begin{equation}\label{eq:fbf}
B^H(\x)=\int_{\R^2}\frac{e^{i \x\cdot \xib}-1}{\|\xib\|^{H+1}}d\widehat W(\xib),\end{equation}
where $d\widehat W$ is a complex Brownian measure and $\x\cdot \xib$ denotes the dot-product on $\R^2$. The Hurst index $H$ is a fundamental parameter of the FBF, as an indicator of the texture roughness. The greater $H$ is, the smoother the resulting texture is, as can be seen in Fig.~\ref{fig:fbf}.

%

\subsection{General anisotropic self-similar Gaussian fields}

In order to introduce anisotropy in this model, Bonami and Estrade \cite{Bonami2003} replaced the Hurst index $H$ in \eqref{eq:fbf}
by a function of the direction of $\xib$ and then derived a new class of {\it Anisotropic Fractional Brownian Field} (AFBF) by: 
\begin{equation}\label{eq:afbf}X(\x)=\int_{\R^2}\frac{e^{i \x\cdot \xib}-1}{\|\xib\|^{h(\arg \xib)+1}}d\widehat W(\xib).\end{equation}
More generally, a larger class of anisotropic models can be defined as
\begin{equation}\label{eq:afbfbis}
X(\x)=\int_{\R^2}(e^{i \x\cdot \xib}-1)f^{1/2}(\xib)~d\widehat W(\xib),
\end{equation}
where the spectral density $f$ is of the form
\begin{equation}\label{eq:density}
f^{1/2}(\xib)=c(\arg \xib)\|\xib\|^{-h(\arg \xib)-1}.
\end{equation}
Here, $c$ and $h$ are two $\pi$-periodic functions, defined on the interval $(-\pi/2,\pi/2]$ with ranges satisfying $c((-\pi/2,\pi/2])\subset \R^+$ and $h((-\pi/2,\pi/2])\subset (0,1)$.
When $c$ and $h$ are both constant, we recover a FBF of order $H\equiv h$.

To define {\it stationary anisotropic models} with \textit{global} orientation $\alpha_0$, one can set  $h\equiv H$ in~(\ref{eq:density}) and:
 \begin{equation}\label{eq:elem}
 c_{\alpha_0,\alpha}(\arg(\xib))=\mathds{1}_{[-\alpha,\alpha]}(\arg(\xib)-\alpha_0),
 \end{equation}
 for some $0<\alpha\leqslant \pi/2$.
 Note that we then recover the \textit{elementary fields} of \cite{Bierme2012}, which are a particular case of AFBF. When $\alpha=\pi/2$, this model corresponds to the usual isotropic FBF of Hurst index $H$ (Fig.~\ref{fig:fbf}), but as soon as $0<\alpha<\pi/2$, the field is no longer isotropic, since the non-zero frequency arguments are restricted between $-\alpha+\alpha_0$ and $\alpha+\alpha_0$. 
 
\begin{figure}[t]
\begin{minipage}[b]{.48\linewidth}
  \centering
  \centerline{\includegraphics[width=4.0cm]{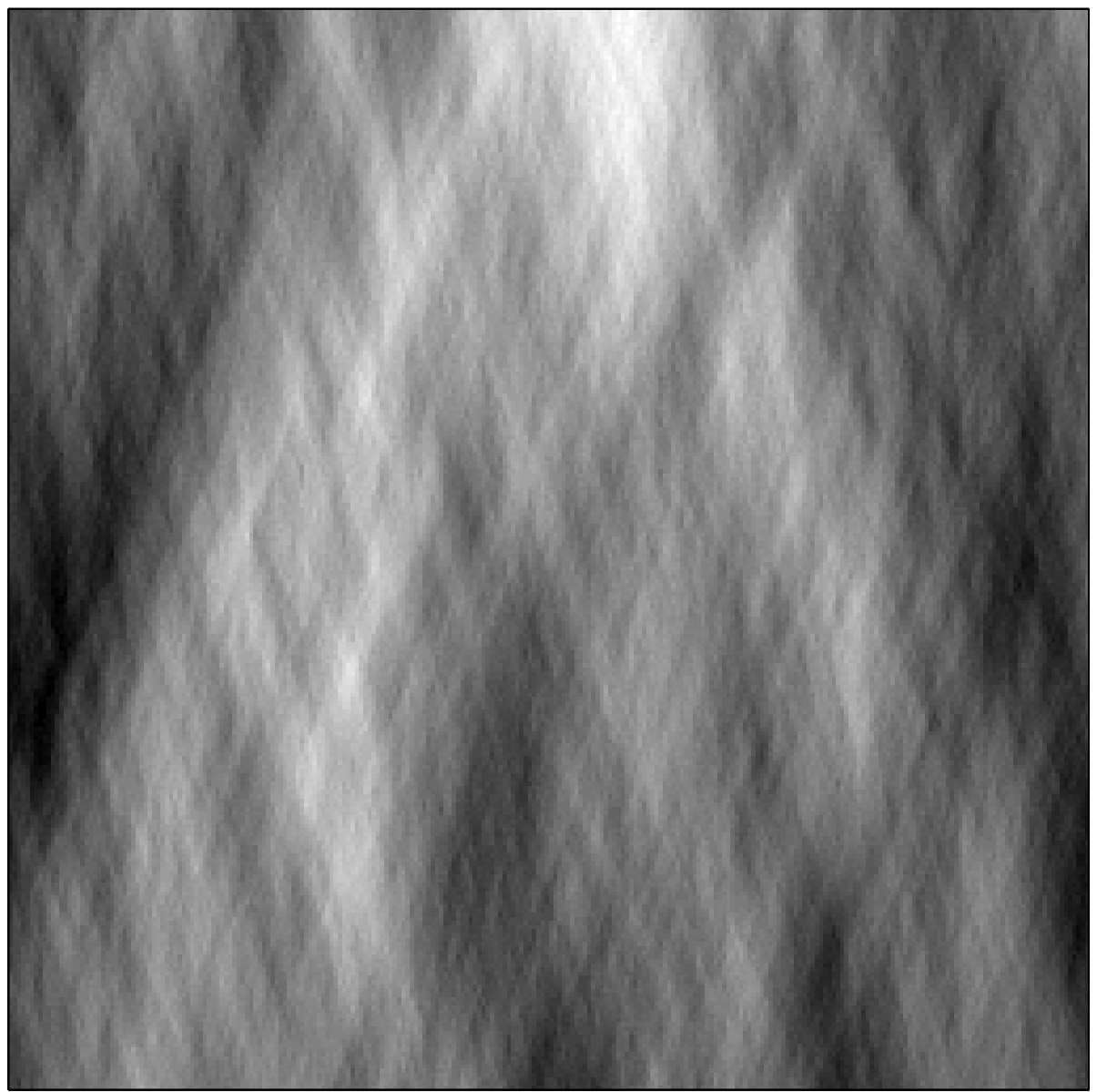}}
  \centerline{(a)}\medskip
\end{minipage}
\hfill
\begin{minipage}[b]{0.48\linewidth}
  \centering
  \centerline{\includegraphics[width=4.0cm]{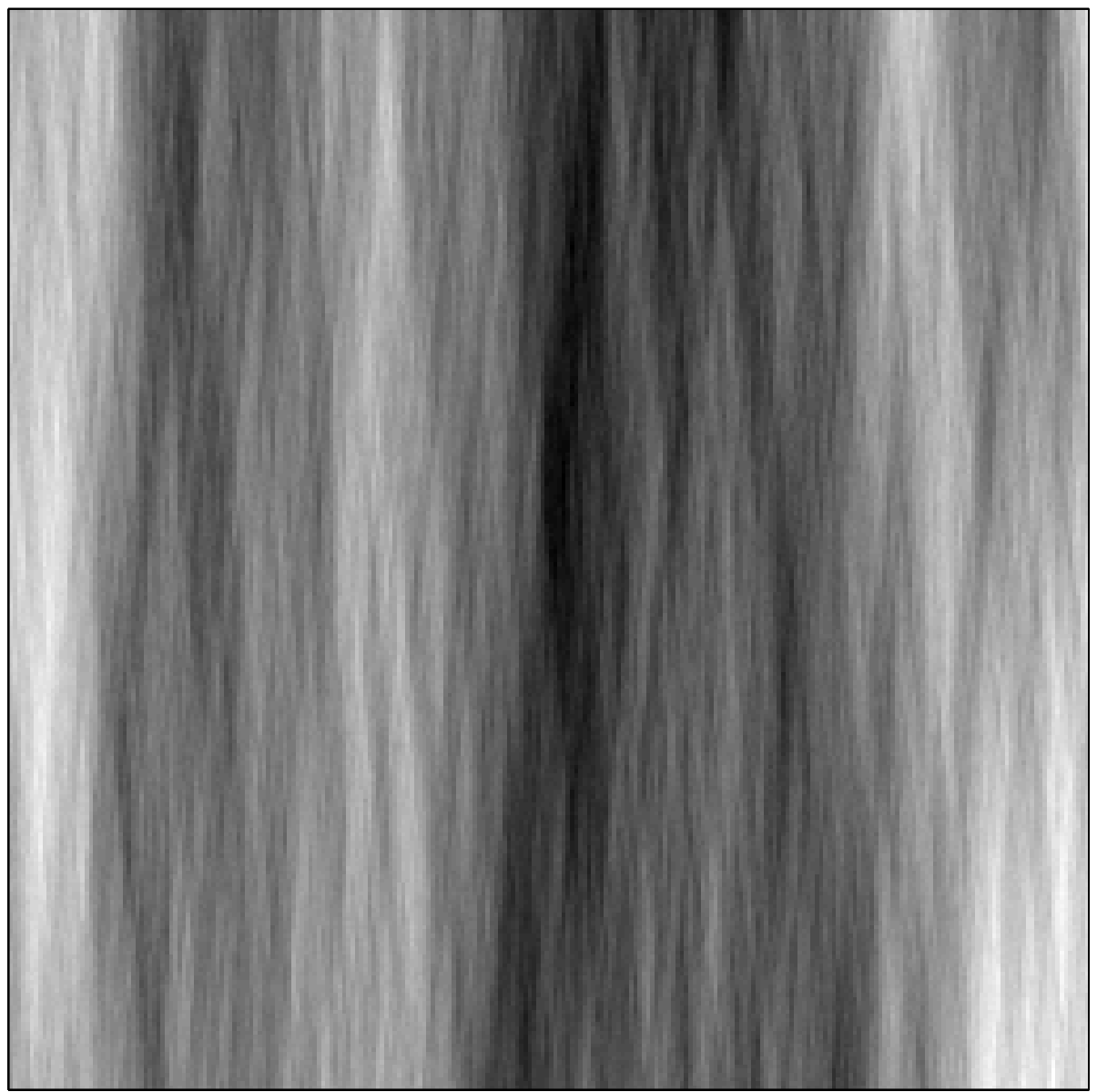}}
  \centerline{(b)}\medskip
\end{minipage}
\caption{AFBF 
with $H=0.5$, $\alpha_0=0$ and (a) $\alpha=\pi/6$, (b) $\alpha=\pi/24$.}
\label{fig:elem}
\end{figure}

Simulation algorithms for Gaussian fields use the covariance function \cite{brouste2007}. But their high complexity is a real problem to produce large textures, and the covariance function is not explicitly known in general case. With respect to the AFBF, a recent fast method has been proposed in \cite{Bierme2012}, called the turning band method, and used here to simulate the textures of Fig.~\ref{fig:elem}, with global orientation $\alpha_0=0$.
Remark that the more the sector $\alpha$ decreases to $0$, the more the frequencies concentrate along the horizontal axis, so the resulting texture appears vertically oriented, as a consequence of the Fourier transform properties. For small $\alpha$, we obtain a strongly stationary oriented texture in the direction orthogonal to $\alpha_0=0$ like in Fig.~\ref{fig:elem} (b).

\begin{figure*}[t]
  \centering
  \includegraphics[scale=0.7]{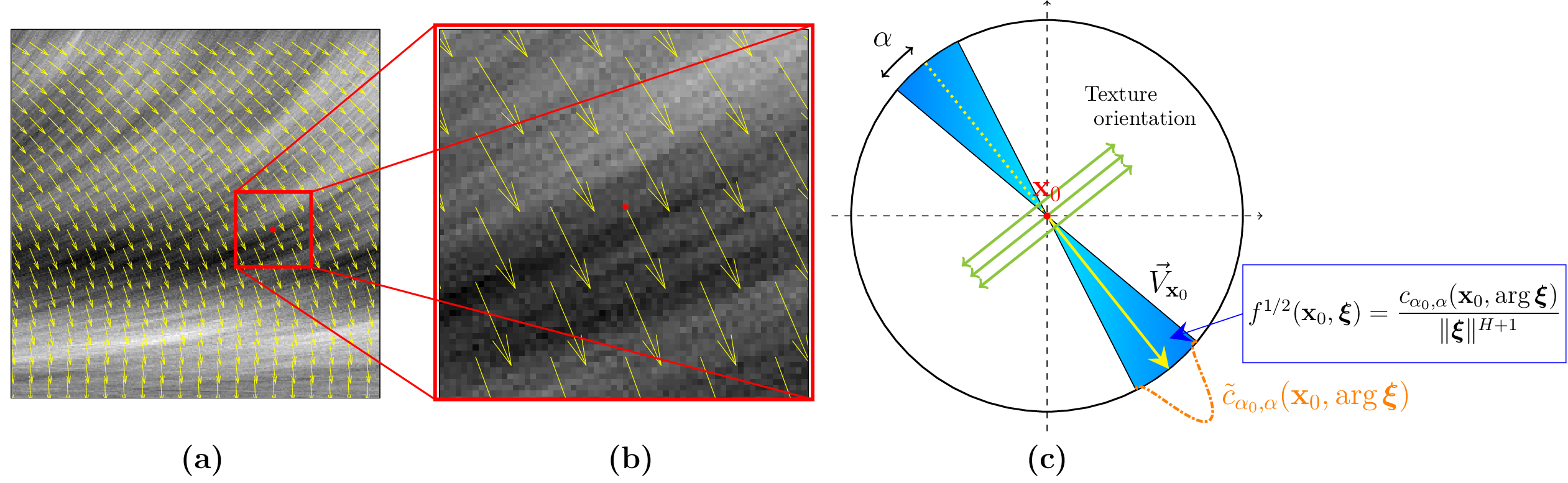}
  \caption{(a) Texture resulting from the vector field orientation $\vec{V}_{(x,y)}^1$ in yellow, (b) zoom around the red point $\x_0=(x,y)$ that shows locally an oriented elementary field, (c) diagram illustrating each parameter of the LAFBF model.}
  \label{fig:summary}
\end{figure*}


\section{A new class of Gaussian field with prescribed orientation}
\label{sec:LAFBF}

\subsection{Definition}

We now define our new Gaussian model as a local version of the elementary field defined in~(\ref{eq:afbfbis}) with density given by~(\ref{eq:density}) with $h\equiv H$ and $c$ as in (\ref{eq:elem}). More precisely, we define the following Gaussian field, that we call {\it Locally Anisotropic Fractional Brownian Field} (LAFBF):
\begin{equation}\label{eq:LAFBF}
X(\x)=\int_{\R^2}(e^{i\x\cdot\xib}-1)f^{1/2}(\x,\xib)~d\widehat{W}(\xib),
\end{equation}
with \begin{equation}\label{eq:dens}f^{1/2}(\x,\xib)=c_{\alpha_0,\alpha}(\x,\arg \xib)\|\xib\|^{-H-1},\end{equation}
\begin{equation}\label{eq:cfunc}c_{\alpha_0,\alpha}(\x,\arg \xib)=\mathds{1}_{[-\alpha,\alpha]}(\arg(\xib)-\alpha_0(\x)),\end{equation} $\alpha_0$ being now a differentiable function on $\R^2$.
Our Gaussian field is derived from AFBF in a similar way than \textit{Multifractional Brownian Motion} (MBM) from FBM, since we replace the orientation parameter $\alpha_0$ with a function depending on the spatial location $\x$, whereas in the case of MBM the Hurst index $H$ of FBM was allowed to vary spatially \cite{peltier1995}.


\subsection{Tangent fields at every point}

To describe the local properties and simulate our new class of Gaussian fields, we shall use the notion of tangent fields that we now briefly review. Recall that the random field $X$ is locally asymptotically self-similar of order $H\in (0,1)$ if for any $\h\in\R^2$ the random field
\[
\frac{X(\x_0+\rho \h)-X(\x_0)}{\rho^H},
\]
admits a non-trivial limit in law $Y_{\x_0}$ as $\rho\to 0$ (see \cite{benassi1997}, and \cite{falconer2002,falconer2003} for a more general definition). The field $Y_{\x_0}$ is then called the tangent field of $X$ at $\x_0$. Roughly speaking, the random field $X$ admits the tangent field $Y_{\x_0}$ at a given point $\x_0$ if it \textbf{behaves locally as} $Y_{\x_0}$ when $\x\to \x_0$. This notion has been first introduced in \cite{benassi1997} to describe the local behavior of Multifractional Brownian Motion (which behaves locally as a FBM).

We can prove that the LAFBF $X$ of \eqref{eq:LAFBF} admits a tangent field $Y_{\x_0}$ at any point $\x_0\in \R^2$ defined as:
\begin{equation}\label{eq:tangent}
Y_{\x_0}(\x)=\displaystyle \int_{\R^2}(e^{i \x\xib}-1)f^{1/2}(\x_0,\xib)~d\widehat{W}(\xib)\;.
\end{equation}
We observe that the tangent field $Y_{\x_0}$ is no more and no less than an elementary field using the terminology of \cite{Bierme2012}. This result shall be crucial when simulating this Gaussian model as detailed in the next section.


\subsection{Simulation of Locally Anisotropic Fractional Brownian Field}

\textbf{Simulation of tangent fields}. The simulation of a LAFBF will first require the simulation of a tangent field at every point $\x_0$. We follow below the methodology of \cite{Bierme2012} using the turning bands method.

\textit{-- Discrete formulation of the tangent field}\\
By a change of variable in polar coordinates, one can derive an integral expression for the variogram of $Y_{\x_0}$: 
\begin{equation}\label{eq:vario}
\begin{array}{rcl}v_{Y_{\x_0}}(\x)&=&\displaystyle \frac{1}{2}\int_{\R^2}|e^{i \x \cdot \xib}-1|^2f(\x_0,\xib)d\xib \\ &=& \displaystyle \frac{1}{2}\gamma(H)\int_{-\pi/2}^{\pi/2}c_{\alpha_0,\alpha}(\x_0,\theta)~|\x\cdot \u(\theta)|^{2H}d\theta\end{array},
\end{equation}
where $\u(\theta)=(\cos \theta,\sin \theta)$ and $\gamma(H)=\frac{\pi}{H\Gamma(2H)\sin(H\pi)}$. \\
The integral (\ref{eq:vario}) is of the form $\int_{-\pi/2}^{\pi/2}\tilde v_{\theta}(\x\cdot \u(\theta))d\theta$ with\\
$\tilde v_{\theta}=\frac{1}{2}\gamma(H)c_{\alpha_0,\alpha}(\x_0,\theta)|\cdot |^{2H}$. Ignoring the factor \\
$\frac{1}{2}\gamma(H)c_{\alpha_0,\alpha}(\x_0,\theta)$, we recognize that $\tilde v_{\theta}$ is the variogram of a FBM of order $H$.  Consequently, $Y_{\x_0}$ can be viewed as a sum of independent FBM rotating around the origin. Discretizing $\theta$ in an ordered set $(\theta_i)_{1\leqslant i\leqslant n}$ of $n$ band orientations, and let be $(\lambda_i)_{1\leqslant i\leqslant n}$  the associated band weights $\lambda_i=\theta_{i+1}-\theta_i$, the {\it turning band fields} take the form
\vspace{-0.2cm}
 \begin{equation}\label{eq:turning}Y_{\x_0}^{[n]}(\x)=\gamma(H)^{\frac{1}{2}}\sum_{i=1}^n\sqrt{\lambda_ic_{\alpha_0,\alpha}(\x_0,\theta_i)}B_i^{H}(\x\cdot \u(\theta_i)),
 \end{equation} 
 where the $B_i^{H}$'s are $n$ independent FBM of order $H$. This discrete version is a good approximation, provided $\displaystyle \max_{i}\lambda_i \leqslant \varepsilon$ for $\varepsilon$ sufficiently small.
 
 \begin{figure}[t]

  \centering
  \includegraphics[scale=0.5]{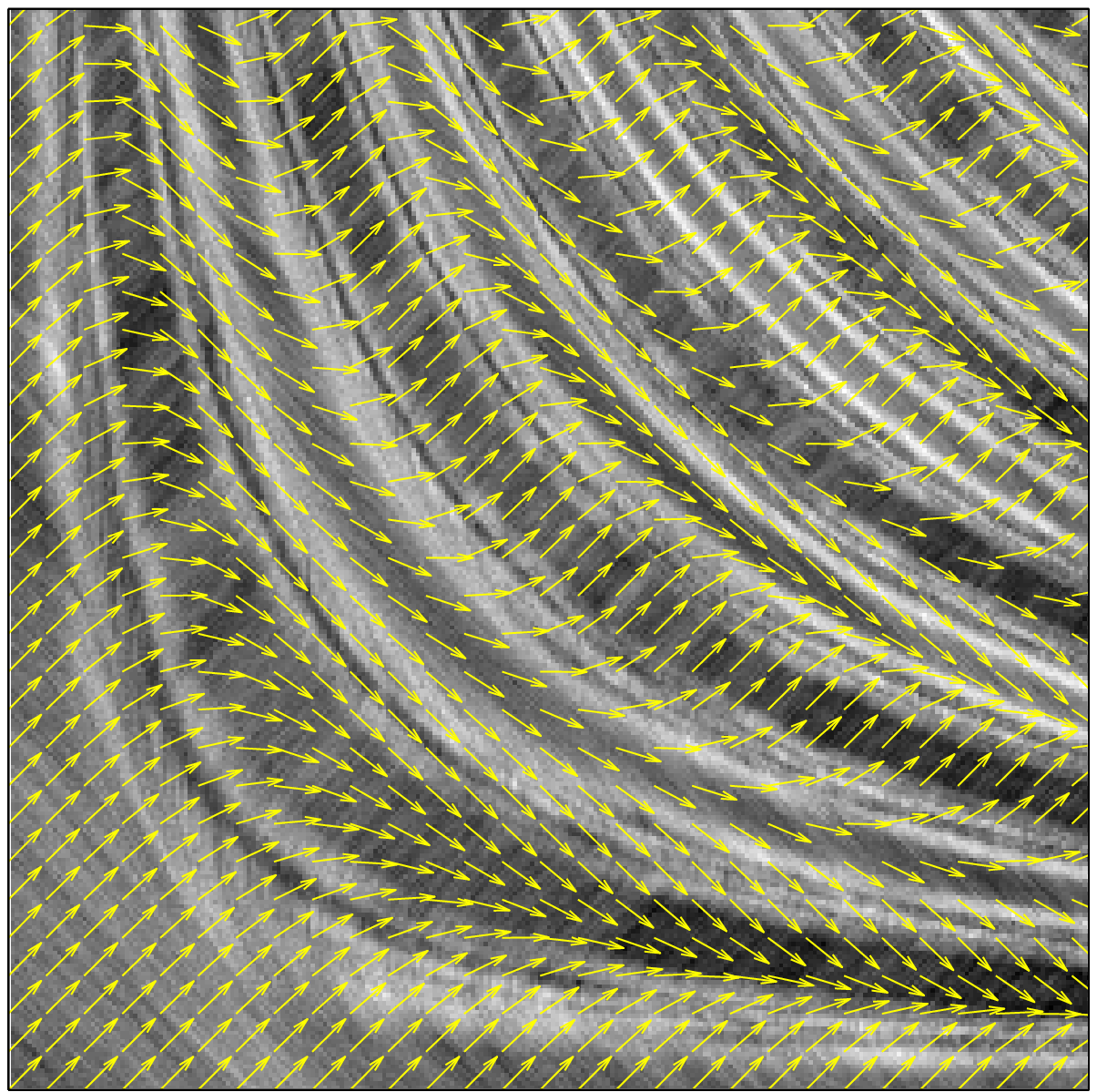}
\caption{Texture resulting from the vector field $\vec{\V}_{(x,y)}^2$.}
\label{fig:vague}
\end{figure}

\textit{-- Simulation along particular bands}\\
 In practice, we consider a discrete grid $r^{-1}\mathbb{Z}^2\cap [0,1]^2$ with $r=2^k-1, k\in \mathbb{N}^{\star}$.  Following \cite{Bierme2012}, we choose $\theta_i$ such that $\tan(\theta_i)=\frac{p_i}{q_i}$, with $p_i,q_i\in \mathbb{Z}$, then $B_i^H(\x\cdot \u(\theta_i))$ becomes \begin{equation}\label{eq:trick}\begin{array}{c}\displaystyle \left\{B_i^{H}\left(\frac{k_1}{r}\cos \theta_i+\frac{k_2}{r}\sin \theta_i\right);0\leqslant k_1,k_2\leqslant r\right\}\overset{\mathcal{L}}{=}\\ \displaystyle
 \left(\frac{\cos \theta_i}{rq_i}\right)^{H}\{B_i^{H}(k_1q_i+k_2p_i);0\leqslant k_1,k_2\leqslant r\}\end{array},\end{equation} and then can be generated using the fast algorithm of Perrin \textit{et al.} \cite{perrin2002} on a regular grid.

\textit{-- Dynamic choice of discrete bands}\\
Finally, the choice of the bands orientations $(\theta_i)_{1\leqslant i\leqslant n}$ is governed by the global computational cost of the $B_i^H$, within dynamic programming  \cite{Bierme2012}.\medskip

\begin{algorithm}[b]
 \caption{Simulation of the LAFBF}
 \begin{algorithmic}[1]
 \INPUT $r=2^k-1$, $H$, $\alpha_0$, $\alpha$, $\epsilon$
 \OUTPUT $X$ LAFBF of size $(r+1)\times (r+1)$
   \STATE $(p_i,q_i)_{1\leqslant i\leqslant n} \leftarrow \text{DynamicBandsChoice}(r,\varepsilon)$
   \STATE Compute and sort angles $(\theta_i)_{1\leqslant i\leqslant n}$ : $\theta_i\leftarrow \text{atan2}(p_i,q_i)$
   \STATE Compute width bands $(\lambda_i)_{1\leqslant i\leqslant n}$ : $\lambda_i\leftarrow \theta_{i+1}-\theta_i$
   \STATE Generate $n$ FBM : $B_i^H\leftarrow \text{circFBM}(r(|p_i|+|q_i|),H)$
   \STATE Initialization : $X\leftarrow 0$
   \FOR {\textbf{all} $(k_1,k_2)$}
		\FOR{$i=1$ \TO $n$}
		        \STATE $\omega_i\leftarrow \sqrt{\lambda_i\gamma(H)c_{\alpha_0,\alpha}((k_1,k_2),\theta_i)}\left(\frac{\cos \theta_i}{rq_i}\right)^H$
   			\STATE $X(k_1,k_2)\leftarrow X(k_1,k_2)+\omega_iB_i^H(k_1q_i+k_2p_i)$
		 \ENDFOR
   \ENDFOR
 \end{algorithmic}
 \end{algorithm}

\textbf{Simulation of the LAFBF}. As observed in  \cite{peltier1995} for the MBM, a Gaussian field can be simulated from its tangent fields. The LAFBF behaving locally like its tangent fields, for every pixel $\x_0$, we assign $X(\x_0)=Y_{\x_0}^{[n]}(\x_0)$. The pseudocode of the algorithm is given below, and the corresponding Matlab code is available on the webpage \cite{polisanoMatlab}. A preprocessing step (instructions 1,2,3,4 in the pseudocode), which does not depend on the expected local orientations, includes the dynamic choice and sorting of discrete bands, and the simulation of the $n$ FBM. These steps are executed once and for all. The rest of the algorithm is of complexity $O(r^2\log n)$. Indeed, at each point $(k_1,k_2)$, a turning band $\theta_i$ contributes to $X(k_1,k_2)$ if and only if $c_{\alpha_0,\alpha}((k_1,k_2),\theta_i)\neq 0$, i.e $|\theta_i-\alpha_0((k_1,k_2))|\leqslant \alpha$. Thus, since the array $\theta_i$ is sorted, one such index $i$ is founded using a binary search, and then the others in its neighborhood.


\begin{figure}[t]
\centering
\includegraphics[scale=0.5]{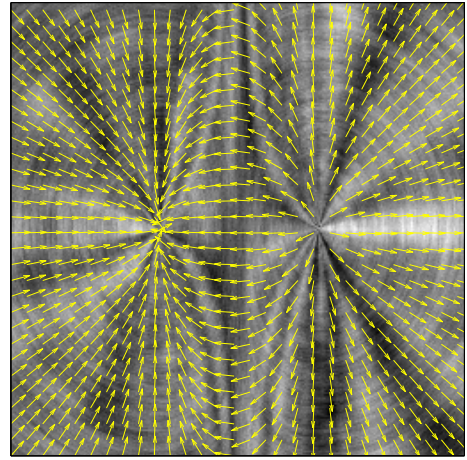}
\caption{Texture resulting from the vector field $\vec{\V}_{(x,y)}^3$.}
\label{fig:magnet}
\end{figure}


\subsection{Oriented texture synthesis}

The parameters used in simulations are $r=255$, $H=0.2$, $\alpha=10^{-1}$, and $\varepsilon=10^{-2}$. To avoid numerical artifacts due to the discrete formula (\ref{eq:turning}) we consider a regularized version $\tilde c_{\alpha_0,\alpha}$ of the indicator function  $c_{\alpha_0,\alpha}$, typically a Gaussian. For $\alpha_0$ constant, we recover the results of \cite{Bierme2012} (see Fig.~\ref{fig:elem}). We present now realizations of textures with prescribed local orientation at each point $\x_0$, given by a vector field $\vec{\V}_{\x_0}=\u(\alpha_0(\x_0))$. Fig.~\ref{fig:summary}(a) displays a texture resulting from the vector field $\vec{\V}_{(x,y)}^1=(\cos(-\pi/2+y),\sin(-\pi/2+y))$. A zoom around a point $\x_0$ (in red, Fig.~\ref{fig:summary}(b)) shows that locally a LAFBF behaves as an elementary field. Fig.~\ref{fig:summary}(c) sketchs the local density function at $\x_0$ and the different parameters. We then consider two others types of vector fields, $\vec{\V}_{(x,y)}^2=(\cos(\cos(36xy)),\sin(\cos(
36xy)))$ and $\vec{\V}_{(x,y)}^3=\nabla F(x,y)$ with $F(x,y)=(4x-2)e^{-(4x-2)^2-(4y-2)^2}$, with the resulting textures in Figs.~\ref{fig:vague},\ref{fig:magnet}. As expected, the textures obtained with our approach present local anisotropic behavior, with a direction orthogonal to the vector field. Moreover, the simulation of a $256\times 256$ texture takes only a few seconds.


\section{Conclusion}

We introduced a new stochastic model defined as a local version of an anisotropic fractional Brownian field. We took advantage of tangent field formulation and the turning bands method to provide an efficient algorithm to simulate textures with prescribed local orientations. We are currently improving the method by removing numerical artifacts which appear for greater values of the Hurst index H. Future extensions of our model include Gaussian fields whose Hurst index and local orientation may vary spatially. A forthcoming work will focus on its application to natural texture characterization and classification, as well as cartoon-texture image decomposition \cite{ono2014cartoon}.

%
%


\bibliographystyle{IEEEbib}
\bibliography{refs}

\end{document}